\documentclass[final,12pt]{elsarticle}
\usepackage{xcolor,soul,cancel}
\usepackage[color]{showkeys}
\usepackage{amssymb}
\usepackage{amsmath}
\usepackage{amsthm}
\usepackage{graphicx}
\usepackage{float}
\journal{arXiv}

\definecolor{refkey}{rgb}{0,1,1}
\definecolor{labelkey}{rgb}{1,0,0}

\numberwithin{equation}{section}

\newtheorem{thm}{Theorem}
\newtheorem{lem}{Lemma}

\newtheorem{prop}[thm]{Proposition}

\newcommand{\abs}[1]{\left\vert#1\right\vert}
\newcommand{\rk}{\operatorname{rank}}

\newcommand{\re}{\operatorname{Re}} 
\newcommand{\real}{\Re} 
\newcommand{\imag}{\Im} 
\newcommand{\conj}{\overline}

\begin{document}

\begin{frontmatter}

\title{{The Gau-Wang-Wu conjecture on partial isometries holds in the 5-by-5 case}
\tnoteref{support}}

\author[nyuad]{Ilya M. Spitkovsky}
\ead{ims2@nyu.edu, ilya@math.wm.edu, imspitkovsky@gmail.com}
\address[nyuad]{Division of Science and Mathematics,
New York  University Abu Dhabi (NYUAD), Saadiyat Island,
P.O. Box 129188 Abu Dhabi, United Arab Emirates}
\author[nyuad]{Ibrahim~Suleiman}
\ead{is1647@nyu.edu}
\author[freu]{Elias~Wegert}
\ead{wegert@math.tu-freiberg.de}
\address[freu]{Institute of Applied Analysis, TU Bergakademie Freiberg,
09596 Freiberg, Germany}
\tnotetext[support]{The results are partially based on the Capstone
project of [IS] under the supervision of [IMS]. The latter was also
supported in part by Faculty Research funding from the Division of Science
and Mathematics, New York University Abu Dhabi. }
\begin{abstract}
Gau, Wang and Wu in their LAMA'2016 paper conjectured (and proved
for $n\leq 4$) that an $n$-by-$n$ partial isometry cannot have a circular
numerical range with a non-zero center. We prove that this statement holds also
for $n=5$.
\end{abstract}

\end{frontmatter}

\section{Introduction}
Let $A$ be a bounded linear operator acting on a Hilbert space ${\mathcal H}$.
The {\em numerical range} (a.k.a. the {\em field of values}, or the {\em
Hausdorff set}) of $A$ is defined as
\[
W(A):= \{\langle{Ax,x}\rangle\colon x\in\mathcal H,\ \left\Vert x\right\Vert=1 \},
\]
where $\langle{.,.}\rangle$ stands for the scalar product on $\mathcal H$, and
$\Vert{.}\Vert$ is the respective norm. Clearly, $W(A)$ is a subset of 
$\{z\in\mathbb{C} \colon \abs{z}\leq\left\Vert A\right\Vert\}$, 
closed in the case
$\dim\mathcal H<\infty$ (but not necessarily in general) and invariant under
unitary similarities of $A$.
By the celebrated Toeplitz-Hausdorff theorem, $W(A)$
is always convex --- see e.g. \cite{GusRa} or \cite[Chapter 1]{HJ2} and the
references therein for a comprehensive treatment of the subject.

An operator $A$ is a {\em partial isometry} if it preserves norms of
the vectors from the orthogonal complement of its kernel:
\[
\left\Vert{Ax}\right\Vert=\left\Vert{x}\right\Vert \text{ whenever } x\in \ker (A)^\perp.
\]
It was conjectured by Gau, Wang, and Wu \cite{GWW} that

{\sl If $A$ is a partial isometry acting on a finite-dimensional space
$\mathcal H$ and $W(A)$ is a circular disk, then the latter is necessarily
centered at the origin.}

In what follows, we will refer to this statement as (GWW). According to
\cite[Theorem 2.2]{GWW}, (GWW) holds for $\dim\mathcal H:=n\leq 4$. To the best
of our knowledge, the question remains open for higher values of $n$. Very
recently, the case $n=5$ was treated in \cite{NaBe} under some additional
restrictions on $A$. In our paper we show that (GWW) for $n=5$ holds
unconditionally.
\begin{thm}\label{th:main}
Let $A$ be a partial isometry acting on a 5-dimensional space,
and assume that its numerical range $W(A)$  is a circular
disk. Then this disk is centered at the origin.
\end{thm}

\section{General remarks}
Denoting by $\mathcal R(A)$ the range of an operator $A$, observe that the
closure $\mathcal L$ of $\mathcal R(A)+\mathcal R(A^*)$ is a reducing subspace
of $A$. Let $A_0$ stand for the restriction of $A$ onto $\mathcal L$; the
restriction of $A$  onto $\mathcal L^\perp=\ker(A)\cap\, \ker(A^*)$ is the zero
operator.
Note that $\mathcal R(A)+\mathcal R(A^*)$ is
automatically closed if $A$ has finite rank  --- the only case we are
interested in. If $A$ is a partial isometry, then $A_0$ is a partial isometry
as well. If, in addition, $W(A)$ is a circular disk, then so is $W(A_0)$. Since
$\dim\mathcal L\leq 2\rk A$, \cite[Theorem 2.2]{GWW} immediately implies the
following

\begin{prop} \label{th:r1}
{\em (GWW)} holds for partial isometries of rank one and two, independent of
the dimension of $\mathcal H$.
\end{prop}

Note also that (GWW) is tautologically correct when $\dim\mathcal L=\rk A$ is
finite. Indeed, in this case $A_0$ is unitary, and so $W(A_0)$ and $W(A)$ are
polygons, never circular disks. When $n=5$, we therefore only need to consider
two remaining values of $\rk A$, namely 3 and 4. This is done in
Sections~\ref{s:3} and~\ref{s:4}, respectively.

Moving forward, we will not distinguish between an operator $A$ acting on a
finite-dimensional space and its matrix representation with respect to an
orthonormal basis. The invariance of $W(A)$ under unitary similarities is the
key to choosing the simplest matrix form possible.

The other ingredient of the reasoning involves the notion of the so called
\emph{\em numerical range generating} (or \emph{Kippenhahn}) \emph{polynomial}
$P_A(\theta,\lambda)$ of $A$. This, by definition, is the \emph{characteristic
polynomial} of $\re (e^{i\theta}A)$, where, as usual, $\re X$ stands for the
hermitian part of $X$. Denoting by $\lambda_j(\theta)$ the roots of this
polynomial, recall that
\[
\{e^{-i\theta}(\lambda_j(\theta)+i\mathbb{R})\colon\ j=1,\ldots,n;\ \theta\in (-\pi,\pi]\}
\]
is the family of tangent lines to a certain algebraic curve $C(A)$, the convex
hull of which is $W(A)$, as was for the first time observed in \cite{Ki} (see
also its more accessible translation \cite{Ki08}). Respectively, $C(A)$ is
called the {\em numerical range generating curve}, or the \emph{Kippenhahn
curve}, of $A$.

We state the pertinent (well known, and also easily verifiable) result for
convenience of reference.
\begin{lem} \label{th:cir}
The Kippenhahn curve $C(X)$ of a given square matrix $X$ contains the circle of
radius $r$ centered at the origin if and only if
$P_X(\theta,\lambda)$ is divisible by $\lambda^2-r^2$ for
all $\theta$.\end{lem}
In Section~\ref{s:3}, assuming that $n=5$ and $\rk A=3$,
we will show that there are no such $a\neq 0$ that make $P_{A-aI}$ divisible by
$\lambda^2-r^2$ for some $r>0$, thus proving (GWW) in that case.
 For $\rk A=4$ the situation is slightly different, since then such
$a$ and $r$ do exist. Instead, a general description of such matrices is given,
and then it is observed that their numerical ranges are not circular
nevertheless.
This is happening since, for such matrices, $C(A)$ does not
lie completely inside the circle $C=\{ z \colon \abs{z-a}=r\}$.

\section{Rank three partial isometries} \label{s:3}
\begin{prop}\label{th:n5k3}
Let $A$ be a rank $3$ partial isometry acting on a $5$-dimensional space. Then
$C(A)$ cannot contain a circle with a non-zero center.
\end{prop}

\begin{proof}
We can choose an orthonormal basis in such a way that the respective matrix of
$A$ (also denoted by $A$) has zeros in the first two columns, and the last
three columns forming an orthonormal set. Invoke the Schur lemma to put the
lower right 3-by-3 block of $A$ in an upper triangular form. Then use a unitary
similarity of the upper left 2-by-2 block to set the (2,3) entry of $A$ at
zero.

Suppose now that $C(A)\supset\{z\colon \abs{z-a}=r\}$ for some $a\in\mathbb{C} , r>0$.
According to \cite{Wu11}, $a$ is an eigenvalue of $A$ of algebraic multiplicity
exceeding its geometric multiplicity. If $a\neq 0$, multiplying $A$ by a
unimodular scalar (under which it remains a partial isometry),
we may assume without loss of generality that $0<a\le 1$.
Furthermore, yet another (this time, diagonal) unitary similarity can be used
to make the entries $a_{13},a_{14},a_{15}$ and $a_{24}$ of $A$ non-negative.
As a result,
\[
A=\begin{bmatrix} 0 & 0 & \sqrt{1-a^2} & ta & sa \\
0 & 0 & 0 & u & v \\ 0 & 0 & a & -t\sqrt{1-a^2} & -s\sqrt{1-a^2} \\
0 & 0 & 0 & a & w \\ 0 & 0 & 0 & 0 & z \end{bmatrix},
\]
where $s,t,u=\sqrt{1-a^2-t^2}\geq 0$, $s^2+\abs{v}^2+\abs{w}^2+\abs{z}^2=1$ and
\begin{equation}\label{st}
st+uv+aw=0.
\end{equation}
Let now
\[
X:=A-aI=\begin{bmatrix} -a & 0 & \sqrt{1-a^2} & ta & sa \\
0 & -a & 0 & u & v \\ 0 & 0 & 0 & -t\sqrt{1-a^2} & -s\sqrt{1-a^2} \\
0 & 0 & 0 & 0 & w \\ 0 & 0 & 0 & 0 & z-a \end{bmatrix},
\]
and write the Kippenhahn polynomial of $X$ as
\begin{equation}\label{Kp3}
P_X(\theta,\lambda)=-\lambda^5+\sum_{j=0}^4A_j(\theta)\lambda^j.
\end{equation}
The explicit formulas for $A_j$, obtained with the use of SageMath, are rather
cumbersome and thus relegated to the Appendix.

Since $C(X)=C(A)-a$ contains a circle centered at the origin, according to
Lemma~\ref{th:cir}, the radius $r$ of this disk is such that
$P_X(\theta,\pm r)=0$. Equivalently,
\begin{equation}\label{r3}
A_0+A_2r^2+A_4r^4=0\text{ and } A_1+A_3r^2-r^4=0.
\end{equation}
%
With the use of the formulas for $A_j$, the equations \eqref{r3} take the form
\begin{equation}\label{b3}
b_0+b_1\cos\theta+b_2\sin\theta+b_3\cos^3\theta+b_4\cos^2\theta\sin\theta+b_5\cos^2\theta=0
\end{equation}
and
\begin{equation}\label{c3}
c_0+c_1\cos\theta\sin\theta+c_2\cos^2\theta=0,
\end{equation}
where the coefficients $b_j,c_j$ depend on the matrix $A$ and on $r$, but not
on~$\theta$. We will provide the explicit expression only for those of them
used in the reasoning below.

Since the function system
\[
\{1,\,\cos\theta,\,\sin\theta,\,\cos^2\theta,\,\cos^2\theta\sin\theta,\,\cos^3\theta\}
\]
is linearly independent, equations \eqref{b3}, \eqref{c3} can only hold for all
values of~$\theta$ if the coefficients $b_j,c_j$ are equal to zero. But
\[
b_3=r^2a^2(\real z-a),\qquad b_4=-r^2a^2\imag z,
\]
so that $z=a$. Substituting this into the formulas for $A_j$ and
evaluating $c_2$, we get
\begin{align} \label{c32}
c_2 = -a^{2} r^{2} &+ \frac{a}{4}\,\imag(w)\,\big(u\imag(v) + a\imag(w)\big) \\
\nonumber&+ \frac{a}{4}\,\real(w)\,\big(s t + u \real(v) + a \real(w)\big).
\end{align}
Rewriting \eqref{st} as
\[
u\imag v + a\imag w = 0, \qquad st+u\real v + a\real w = 0,
\]
we see that all the summands in the right hand side of \eqref{c32}, except for
the very first one, cancel out. Hence $c_2=-a^2r^2=0$, a contradiction to $a,r>0$.
\end{proof}

\section{Rank four partial isometries} \label{s:4}

\noindent As it happens, a rank four 5-by-5 partial isometry $A$ can have a
circle with a non-zero center as part of its Kippenhahn curve. The location of
this center defines $A$ up to a unitary similarity. The radius of this circle
is therefore defined uniquely. More specifically, the following statement
holds.
\begin{prop} \label{n5k4}
Let $A$ be a rank four 5-by-5 partial isometry. Then $C(A)$ contains a circle
centered at $ae^{i\omega}$ for some $a>0$ if and only if $a<1$ and $A$ is
unitarily similar to
\begin{equation}\label{c45} e^{i\omega}
\begin{bmatrix}0 & \sqrt{1-a^2} & a\sqrt{1-a^2}  & a^2\sqrt{1-a^2}  & a^3\sqrt{1-a^2}  \\
0 & a & a^2-1 & a(a^2-1) & a^2(a^2-1) \\ 0 & 0 & a & a^2-1 & a(a^2-1) \\
0 & 0 & 0 & a & a^2-1\\ 0 & 0 & 0 & 0 & a
\end{bmatrix}.
\end{equation}
If this is the case, then the radius of the circle in question is $r=(1-a^2)/2$.
\end{prop}
\begin{proof}
{\sl Necessity.} Choose an upper triangular matrix representation of $A$, with
the zero first column and the remaining four forming an orthonormal set.
Without loss of generality we may (and, to simplify the notation, will) suppose
that the center of the circle contained in $C(A)$ is positive, i.e.
$e^{i\omega}=1$. At least two of the eigenvalues of $A$ are equal to $a$.
Applying a permutational similarity if needed, we may set $a_{22}=a_{33}=a$.
Also, as in the proof of Proposition~\ref{th:n5k3}, the first row of $A$ may be
arranged to become entry-wise non-negative. The orthonormality of the second
and the third columns of $A$ and their orthogonality to the  fourth and the
fifth then yield
\begin{equation}\label{Aint}
A= \begin{bmatrix}
0 & \sqrt{1-a^2} & a\sqrt{1-a^2}  & xa\sqrt{1-a^2}  & ya\sqrt{1-a^2}  \\
0 & a & a^2-1 & x(a^2-1) & y(a^2-1) \\ 0 & 0 & a & x(a^2-1)/a & y(a^2-1)/a \\
0 & 0 & 0 & u & v \\ 0 & 0 & 0 & 0 & w
\end{bmatrix}
\end{equation}
with $x,y\geq 0$ and some, yet to be determined, parameters $u,v,w\in\mathbb{C} $.

It is time now to invoke Lemma~\ref{th:cir}, according to which the Kippenhahn
polynomial $P_X$ of $X=A-aI$ is divisible by $\lambda^2-r^2$ for some $r>0$.
Writing $P_X$ as
\[
P_X(\theta,\lambda) = -\lambda^5 + \sum_{j=0}^4 B_j(\theta)\,\lambda^j,
\]
we conclude that
\begin{equation}\label{B4}
B_0+B_2r^2+B_4r^4 = 0\quad \text{and}\quad B_1+B_3r^2-r^4=0
\end{equation}
for all values of $\theta$. Using the explicit formulas for $B_j$
(obtained, as in Section~\ref{s:3}, with the use of SageMath, and listed in the
Appendix) observe that
\[
B_0+B_2r^2+B_4r^4=b_0+b_1e^{i\theta}+b_2e^{-i\theta}+b_3e^{3i\theta}+b_4e^{-3i\theta},
\]
\[
B_1+B_3r^2-r^4=c_0+c_1e^{2i\theta}+c_2e^{-2i\theta},
\]
with $b_j,c_j$ expressed in terms of $a,x,y,u,v,w$ and $r$. It follows from
\eqref{B4} that all $b_j,c_j$ are equal to zero. Now we have
\[
b_3=\frac{1}{8}(a^2r^2u+a^2r^2w-a^3r^2-ar^2uw)
\]
and
\begin{equation}\label{c14}
c_1=\frac{1}{8}(-a^{4} - 6a^{2}r^{2} - (a^{2} + 2r^{2} - 1)uw + a^{2} + (a^{3}
+ 4 \, a r^{2} - a)(u + w)).
\end{equation}
Since $b_3$ can be rewritten as $-\frac{1}{8}ar^2(u-a)(w-a)$,
either $u=a$ or $w=a$. Assuming the latter, \eqref{c14} simplifies to
$c_1=\frac{1}{4}ar^2(u-a)$, and so $u=a$ as well.
Similarly, if $u=a$, then $c_1=0$ implies $w=a$. So, in fact $u=w=a$.
Plugging these values of $u,w$ into \eqref{Aint} we obtain
\[
A= \begin{bmatrix}0 & \sqrt{1-a^2} & a\sqrt{1-a^2}  & xa\sqrt{1-a^2}  & ya\sqrt{1-a^2}  \\
0 & a & a^2-1 & x(a^2-1) & y(a^2-1) \\
0 & 0 & a & x(a^2-1)/a & y(a^2-1)/a \\
0 & 0 & 0 & a & v \\
0 & 0 & 0 & 0 & a
\end{bmatrix}.
\]
The unit length of column four implies now that $x=a$. Given that, the
orthogonality of the last two columns allows us to express $v$ in terms of $y$
as $v=y(a^2-1)/a^2$. But the fifth column also has to have length one, so
$y=a^2$ and $v$ simplifies to $a^2-1$. We have thus arrived at \eqref{c45} with
$\omega=0$.

To compute $r$, let us use the expression for $b_1 (=b_2)$ which for the matrix
under consideration takes the form
\[
\frac{1}{32}(a^9-4a^7+6a^5-16ar^4-4a^3+a)
\]
and, up to an inconsequential constant multiple, factors as
\[
(a^4-2a^2+4r^2+1)(a^2+2r-1)(a^2-2r-1)a.
\]
Since $a\leq 1$ and $r>0$, this product vanishes if and only if $r=(1-a^2)/2$.

{\sl Sufficiency.} Direct computations show that, for $A$ given by \eqref{c45}
with $\omega=0$, $0<a<1$ and $r=(1-a^2)/2$, the Kippenhahn polynomial of
$X=A-aI$ factors as
\begin{equation}\label{p45}
-(\lambda^2-r^2)(\lambda^3+\sqrt{1-2r}(\cos\theta)\lambda^2+(r^2-2r)\lambda+r^2\sqrt{1-2r}
\cos\theta).
\end{equation}
\end{proof}

Note that the eigenvalue $e^{i\omega}a$ of the matrix \eqref{c45}
has geometric multiplicity one while the zero eigenvalue is not reducing.
Consequently, this matrix is unitarily irreducible while its Kippenhahn polynomial factors.

An explicit description of $W(A)$ for matrices \eqref{c45} can be obtained
with the use of \eqref{p45}. For our purposes, however, it suffices to observe
that it differs from $D:=\{z\colon \abs{z-a}\leq (1-a^2)/2\}$, and thus is not a
circular disk. To this end, recall that the numerical range of any matrix
contains numerical ranges of its principal submatrices. In particular, $W(A)$
contains the numerical range of its upper left 2-by-2 block. By the Elliptical
Range Theorem (and, once again, setting for simplicity $\omega=0$), the latter
is the elliptical disk with the foci $0$ and $a$ and the minor axis of length $\sqrt{1-a^2}$.
The leftmost point of this ellipse is $(a-1)/2$, and its distance from $a$ is
$(1+a)/2>r$.

This observation, combined with Propositions~\ref{th:r1}--\ref{n5k4}, completes the proof of Theorem~\ref{th:main}. \qed

More can be said about the shape of $W(A)$ of matrices \eqref{c45}. In particular, denoting the second factor of \eqref{p45} by $p(\lambda)$, observe that \[ p(\pm r)=2r^2(\sqrt{1-2r}\cos\theta\pm (r-1)). \]
Consequently, $p(r)<0,\ p(-r)>0$, and so $p$ has roots both smaller than $-r$ and larger than $r$. In other words, $D$ lies in the interior of $W(A)$.

\providecommand{\bysame}{\leavevmode\hbox to3em{\hrulefill}\thinspace}
\providecommand{\MR}{\relax\ifhmode\unskip\space\fi MR }
\providecommand{\MRhref}[2]{%
	\href{http://www.ams.org/mathscinet-getitem?mr=#1}{#2}
}
\providecommand{\href}[2]{#2}

\section*{Appendix}

\noindent
The coefficients of $P_X$ in the rank three case are:
\begin{alignat*}{2}
A_0  &= &&\,a_{01}\,\cos\theta + a_{02}\,\sin\theta, \quad \text{where}\\
&&&a_{01} := \frac{1}{16}\,(a^2-1)\left(a u^2+a\,\vert w\vert ^2
+ u\,\real (v\,\overline{w})-u^2\real(z)\right)\\\
&&&a_{02} := \frac{1}{16}\,(a^2-1)\,\left(u\,\imag (v\,\overline{w})
+ u^2\,\imag(z)\right), \\[1ex]
A_1 &= &&\,a_{10} + a_{11}\cos\theta\sin\theta + a_{12}\cos^2 \theta, \quad \text{where} \\
&&&a_{10} :=  \frac{1}{16} \,\left((a^2-1)(u^2 + \vert v\vert ^2 + \vert w\vert ^2)
        + 2\,s\,t\, u\,\real(v) - t^2\vert v\vert ^2 - s^2u^2\right) \\
&&&a_{11} := \frac{1}{4}\,a\,\left((a^2-t^2-u^2-1)\imag(z)
              + s\,t\,\imag(w)+u\,\imag(w\conj{v})\right)\\
&&&a_{12} := \frac{1}{4}\,a\left((a^2 - t^2 - u^2 - 1)(\real(z)-a)
              + u\,\real(v \conj{w}) + a\,\vert w\vert ^2 + s\,t\,\real(w)\right),\\[1ex]
A_2 &= &&\,a_{21}\cos\theta + a_{22}\sin\theta + a_{23}\cos^2\theta\sin\theta
+ a_{24}\cos^3\theta, \quad \text{where}\\
&&&a_{21} := \frac{1}{4}\,\left(a(s^2+\vert v\vert ^2+2\vert w\vert ^2)
+ (a^2-t^2-u^2-1)(\real z-2a) + st\real w + u\real(v \conj{w}) \right)\\
&&&a_{22} := \frac{1}{4}\,\left((t^2+u^2+1-a^2)\imag(z) - s\,t\,\imag(w)
+ u\,\imag(v\conj w)\right)\\
&&&a_{23} := - a^2\imag(z) \\
&&&a_{24} := a^2(\real(z) - a), \\[1ex]
A_3 &= &&\,a_{30} + a_{31}\cos^2\theta + a_{32}\cos\theta\sin\theta,
\quad \text{where}\\
&&&a_{30} := \frac{1}{4}\,\left(s^2 + t^2 + u^2 + \vert v\vert ^2
+ \vert w\vert ^2 + 1 - a^2\right) \\
&&&a_{31} := a\,(2\,\real(z) - 3\,a)\\
&&&a_{32} := -2\,a\,\imag(z),\\[1ex]
A_4 &=&&\,(\real(z)-3a)\,\cos\theta - \imag(z)\,\sin\theta.
\end{alignat*}

\newpage

\noindent
The coefficients of $P_X$ in the rank four case are:

\begin{alignat*}{2}
B_0  &= &&\,b_{01}\,e^{i\theta} + \conj{b_{01}}\,e^{-i\theta}, \quad \text{where}\\
&&&b_{01} := \frac{1}{32}\,\frac{(1-a^2)^3}{a^2}\,
\left(w x^2 + u y^2 - ax^2 - a y^2 - v x y\right)\\[1ex]
B_1 &= &&\,b_{10} + b_{12}\,e^{2i\theta} + \conj{b_{12}}\,e^{-2i\theta}, \quad
\text{where}\\
&&&b_{10} :=  \frac{1}{4}\,(1-a^2)\left(a^2-a\real u - a\real w
   + \real(w\conj{u})-\frac{\vert v\vert ^2}{2}\right) -\frac{(1-a^2)^3}{16a^2}(x^2+y^2)\\
&&&b_{12} := \frac{1}{8}\,(1 - a^2)(a^2 - a u - a w + u w )\\[1ex]
B_2 &= &&\,b_{21}\,e^{i\theta} + \conj{b_{21}}\,e^{-i\theta}
+ b_{23}\,e^{3i\theta} + \conj{b_{23}}\,e^{-3i\theta}, \quad
\text{where}\\
&&&b_{21} := \frac{1-a^2}{8a^2}(v x y - w x^2 - u y^2)
+ \frac{1-a^2}{a}(4a^2 - 2 a u - 2 a w + x^2 + y^2)\\
&&&\hspace{10mm}+ a^2(u  + w  + 2\,\real u + 2\,\real w - 3a)
+ a\,(\vert v\vert ^2 - 2w\,\real u  - u \conj{w})\\
&&&b_{23} :=\frac{1}{8}(a^2u + a^2w - a u w - a^3) \\[1ex]
B_3 &= &&\,b_{30} + b_{32}\,e^{2i\theta} + \conj{b_{32}}\,e^{-2i\theta}, \quad
\text{where}\\
&&&b_{30} := \frac{1}{4}\left(2-8a^2 - x^2 - y^2 + 4a\,\real(u+w)
- 2\,\real(u\conj w) + \vert v\vert ^2
+ 1 + \frac{x^2}{a^2} + \frac{y^2}{a^2}\right) \\
&&&b_{32} := \frac{1}{4}\,(2 a u + 2 a w - u w - 3a^2)\\[1ex]
B_4 &=&&\,\frac{1}{2}(u+v-3a)\,e^{i\theta}
+ \frac{1}{2}(\conj{u}+\conj{v}-3a)\,e^{-i\theta}.
\end{alignat*}


\begin{thebibliography}{1}
	
	\bibitem{GWW}
	H.-L. Gau, K.-Z. Wang, and P.~Y. Wu, \emph{Circular numerical ranges of partial
		isometries}, Linear Multilinear Algebra \textbf{64} (2016), no.~1, 14--35.
	
	\bibitem{GusRa}
	K.~E. Gustafson and D.~K.~M. Rao, \emph{Numerical range. {T}he field of values
		of linear operators and matrices}, Springer, New York, 1997.
	
	\bibitem{HJ2}
	R.~A. Horn and C.~R. Johnson, \emph{Topics in matrix analysis}, Cambridge
	University Press, Cambridge, 1994, Corrected reprint of the 1991 original.
	
	\bibitem{Ki}
	R.~Kippenhahn, \emph{{\"{U}}ber den {W}ertevorrat einer {M}atrix}, Math. Nachr.
	\textbf{6} (1951), 193--228.
	
	\bibitem{Ki08}
	\bysame, \emph{On the numerical range of a matrix}, Linear Multilinear Algebra
	\textbf{56} (2008), no.~1-2, 185--225, Translated from the German by Paul F.
	Zachlin and Michiel E. Hochstenbach.
	
	\bibitem{NaBe}
	M.~Naimi and M.~Benharrat, \emph{On the circulrar numerical ranges of 5-by-5
		partial isometries}, arXiv.math.FA/2108.04459v1 (2021), 1--17.
	
	\bibitem{Wu11}
	P.~Y. Wu, \emph{Numerical ranges as circular discs}, Appl. Math. Lett.
	\textbf{24} (2011), no.~12, 2115--2117.
	
\end{thebibliography}
\end{document}